\documentclass[10pt]{article}
\usepackage{amsmath,amssymb}
\newtheorem{Th}{Theorem}[section]
\newtheorem{Le}{Lemma}[section]

\newtheorem{Rem}{Remark}[section]

\date{}
\begin{document}

\title{The Dirichlet problem for the nonstationary Stokes system in a polygon}
\author{J\"urgen Rossmann}

\maketitle

\section{Introduction}

Let $\Omega$ be a polygonal domain with the boundary $\partial\Omega=\Gamma\cup {\cal P}$, where
${\cal P}=\{P_1,\ldots,P_n\}$ is a set of corner points and $\Gamma$ is the union of the smooth ($\in C^2$) arcs $\Gamma_1=P_1P_2,\ldots,\Gamma_{n-1}=P_{n-1}P_n,\Gamma_n=P_nP_1$. For simplicity, we assume that $\Omega$
coincides with an angle $K_j$ with opening $\alpha_j$ in a neighborhood of the corner point $P_j$, $j=1,\ldots,n$.
The present paper deals with the initial-boundary value problem
\begin{eqnarray} \label{stokes1}
&&\frac{\partial {\mathfrak u}}{\partial t} - \Delta {\mathfrak u} -\nabla\nabla\cdot{\mathfrak u} + \nabla {\mathfrak p} = {\mathfrak f}, \quad
-\nabla\cdot {\mathfrak u} = {\mathfrak g} \ \mbox{ in } Q_T=\Omega \times (0,T), \\ \label{stokes2}
&& {\mathfrak u}(x,t)=0 \ \mbox{ for }x\in \partial\Omega\backslash {\cal P}, \ 0<t<T, \\ \label{stokes3}
&& {\mathfrak u}(x,0)=0 \ \mbox{ for }x\in \Omega.
\end{eqnarray}
The analogous problem in an infinite angle was already studied in \cite{r-18}. As in \cite{r-18}, the major part of the paper deals with the parameter-depending problem
\begin{eqnarray} \label{par3}
&& s\, u - \Delta u -\nabla\nabla\cdot u+ \nabla p = f, \quad -\nabla\cdot u = g \mbox{ in } \Omega, \\ \label{par4}
&& u=0\ \mbox{ on }\Gamma.
\end{eqnarray}
We are interested in solutions of this problem in the class of the weighted Sobolev spaces $V_\beta^l(\Omega)$
which are defined for nonnegative integer $l$ and real $n$-tuples $\beta=(\beta_1,\ldots,\beta_n)$ as the spaces
(closure of the set $C_0^\infty(\overline{\Omega}\backslash {\cal P})$ with the norm
\[
\| u\|_{V_\beta^l(\Omega)} = \Big( \int_\Omega \sum_{|\alpha|\le l} \prod_{j=1}^n r_j^{2(\beta_j-l+|\alpha|)}\,
  |\partial_x^\alpha u(x)|^2\, dx\Big)^{1/2},
\]
where $r_j(x)$ denotes the distance of the point $x$ from the corner point $P_j$.
We suppose that $f\in V_\beta^0(\Omega)$ and $g\in V_\beta^1(\Omega)$, where the inequalities
\begin{equation} \label{3ct1}
\max\big(1-\mbox{Re}\, \lambda^*(\alpha_j), -\frac \pi{\alpha_j}\big) < \beta_j < \min\big(1,\frac\pi{\alpha_j}\big), \quad \beta_j\not=0
\end{equation}
are satisfied for the components $\beta_j$ of $\beta$. Here $\lambda^*(\alpha)$ is the solution of the equation
\begin{equation} \label{evequation}
\sin(\lambda\alpha) + \lambda\sin\alpha =0
\end{equation}
with smallest positive real part. It is shown in Section 3 (see Theorem \ref{ct1}) that the problem (\ref{par3}), (\ref{par4}) has a unique solution $(u,p)$ satisfying the
condition $\int_\Omega p(x)\, dx =0$ and the estimate
\[
\| u\|_{V_\beta^2(\Omega)} + |s|\, \| u \|_{V_\beta^0(\Omega)} + \| p\|_{W_\beta^1(\Omega)} \le  c\, \Big( \| f\|_{V_\beta^0(\Omega)}
+ \| g\|_{V_\beta^1(\Omega)}  + \ |s|\, \| g\|_{(W_{-\beta}^1(\Omega))^*} \Big)
\]
with a constant $c$ independent of $f,g$ and $s$ provided that $\mbox{Re}\, s \ge 0$ and  $|s|>\gamma_0$, where $\gamma_0$ is a sufficiently large positive number.
Here $W_\beta^1(\Omega)= \{ p\in L_2(\Omega):\ \nabla p \in V_\beta^1(\Omega)\}$.

Applying the inverse Laplace transform, we prove first the existence of solutions of the problem (\ref{stokes1})--(\ref{stokes3}) for $T=\infty$ in weighted
spaces with the additional weight function $e^{-\gamma t}$, $\gamma >\gamma_0$ (see Theorem \ref{ct2}). After this it is not difficult to prove
the solvability of this problem on a finite $t$-interval. The main result of the paper (Theorem \ref{ct3}) ist the following:

{\em Suppose that ${\mathfrak f}\in L_2(0,T;V_\beta^0(\Omega))$, ${\mathfrak g}\in L_2(0,T;V_\beta^1(\Omega))$ and
$\partial_t {\mathfrak g}\in L_2(0,T;(W_{-\beta}^1(\Omega))^*)$, where the inequalities {\em (\ref{3ct1})} are satisfied for all $j$.
Furthermore, we assume that ${\mathfrak g}(x,0)=0$ for $x\in \Omega$ and
\[
\int_\Omega {\mathfrak g}(x,t)\, dx = 0 \ \mbox{ for almost all }t.
\]
Then there exists a uniquely determined solution $({\mathfrak u},{\mathfrak p})$ of the problem {\em (\ref{stokes1})--(\ref{stokes3})} such that
${\mathfrak u}\in L_2(0,T;V_\beta^2(\Omega))$, ${\mathfrak u}_t\in L_2(0,T;V_\beta^0(\Omega))$,
${\mathfrak p} \in L_2(\Omega\times(0,T))$, $\nabla{\mathfrak p}\in L_2(0,T;V_\beta^0(\Omega))$ and}
\[
\int_\Omega {\mathfrak p}(x,t)\, dx =0 \ for \ almost\ all \ t.
\]

\section{The parameter-dependent problem in an angle}

Here, we present a theorem which was proved in \cite[Theorem 4.12 and Corollary 4.13]{r-18}.
Let $K$ be an angle with aperture $\alpha$, i.e.,  $K=\{ x=(x_1,x_2): \ 0<r<\infty,\ 0<\varphi<\alpha\}$, where $r$, $\varphi$ denote the polar
coordinates of the point $x$. The sides of $K$ are the half-lines $\Gamma_1$ and $\Gamma_2$,
where $\varphi=0$ on $\Gamma_1$ and $\varphi=\alpha$ on $\Gamma_2$.
We consider the problem
\begin{eqnarray} \label{par1}
&& s\, u - \Delta u -\nabla\nabla\cdot u+ \nabla p = f, \quad -\nabla\cdot u = g \mbox{ in } K, \\  \label{par2}
&& u=0 \mbox{ on } \Gamma_1\cup \Gamma_2.
\end{eqnarray}
For nonnegative integer $l$ and real $\beta$, we define the weighted Sobolev spaces $V_\beta^l(K)$
as the set of all functions (or vector functions) with finite norm
\begin{equation} \label{Vbeta}
\| u\|_{V_\beta^l(K)} = \Big( \int_K r^{2(\beta-l+|\alpha|)}\, \big| \partial_x^\alpha u(x)\big|^2\, dx\Big)^{1/2}.
\end{equation}
The intersection $V_\beta^2(K) \cap V_\beta^0(K)$ is denoted by $E_\beta^2(K)$. Furthermore, let $W_\beta^1(K)$ be the space with the norm
\[
\| g\|_{W_\beta^1(K)} = \Big( \int\limits_{\substack{K \\ 1<|x|<2}} |g|^2\, dx + \int_K r^{2\beta} |\nabla g|^2\, dx\Big)^{1/2}.
\]

\begin{Th} \label{bt2}
Suppose that $s\not=0$, $\mbox{\em Re}\, s \ge 0$ and
$\max(-\frac \pi\alpha \, ,1-\mbox{\em Re}\, \lambda^*(\alpha))< \beta <\min(2,\frac\pi\alpha)$, $\beta\not=0$,
where $\lambda^*(\alpha)$ is the solution of the equation {\em (\ref{evequation})} with smallest positive real part.
Furthermore, let  $f\in V_\beta^0(K)$ and $g\in V_\beta^1(K)\cap (W_{-\beta}^1(K))^*$. In the case $0<\beta<2$,
we assume that the integral of $g$ over $K$ is zero. Then the the problem
\begin{equation} \label{2bt2}
s\, u - \Delta u -\nabla\nabla\cdot u+ \nabla p = f, \quad -\nabla\cdot u = g \mbox{ in } K, \quad
u=0 \mbox{ on } \Gamma_1\cup \Gamma_2
\end{equation}
has a unique solution $(u,p)\in E_\beta^2(K)\times V_\beta^1(K)$ satisfying the estimate
\begin{equation} \label{1bt2}
\| u\|_{V_\beta^2(K)} + |s|\, \| u\|_{V_\beta^0(K)} + \| p\|_{V_\beta^1(K)}
 \le c\, \Big( \| f\|_{V_\beta^0(K)} + \| g\|_{V_\beta^1(K)}  + |s|\, \| g\|_{(W_{-\beta}^1(K))^*} \Big)
\end{equation}
with a constant $c$ independent of $f$, $g$ and $s$.
\end{Th}

Note that $\mbox{Re}\, \lambda^*(\alpha)> 1/2$ for $\alpha<2\pi$ and  $\mbox{Re}\, \lambda^*(\alpha)> 1$ for $\alpha<\pi$.

\begin{Rem} \label{br1}
{\em In \cite[Theorem 4.12]{r-18}, the assertion of Theorem \ref{bt2} was obtained under the assumption that 
$f\in V_\beta^0(K)$ and $g\in V_\beta^1(K)\cap (V_{-\beta}^1(K))^*$. But this assumption is the same as in Theorem \ref{bt2}
if the integral of $g$ is zero for $0<\beta<2$. 
Indeed, if $\beta<0$, then it follows from Hardy's inequality that the spaces $W_{-\beta}^1(K)$ and $V_{-\beta}^1(K)$ coincide.
Suppose that $\beta>0$, $g\in V_\beta^1(K)\cap (V_{-\beta}^1(K))^*$ and $\int_K g\, dx =0$. Then any $q\in W_{-\beta}^1(K)$ is continuous
at the corner point $x=0$ and
\begin{eqnarray*}
\Big| \int_K g\, q\, dx\Big| & = & \Big| \int_K g\, (q-q(0))\, dx\Big| \le \| g\|_{(V_{-\beta}^1(K))^*} \ \| q-q(0)\|_{V_{-\beta}^1(K)} \\
& \le & c_1\, \| g\|_{(V_{-\beta}^1(K))^*} \ \| q\|_{W_{-\beta}^1(K)}
\end{eqnarray*}
(see \cite[Lemma 7.1.3]{kmr-97}). This means that $g\in (W_{-\beta}^1(K))^*$,
\[
\| g\|_{(W_{-\beta}^1(K))^*} \le c_1\, \| g\|_{(V_{-\beta}^1(K))^*} \le c_2\, \| g\|_{(W_{-\beta}^1(K))^*}  \, ,
\]
and the norm of $g$ in $(W_{-\beta}^1(K))^*$ in the estimate (\ref{1bt2}) can be replaced by the norm in the space $(V_{-\beta}^1(K))^*$.}
\end{Rem}

\begin{Rem} \label{br2}
{\em Suppose that $f\in V_\beta^0(K) \cap V_\gamma^0(K)$, $g\in V_\beta^1(K)\cap V_\gamma^1(K)$ and $g\in \cap (W_{-\beta}^1(K))^* \cap (W_{-\gamma}^1(K))^*$,
where
\[
\max\Big(-\frac \pi\alpha \, ,1-\mbox{Re}\, \lambda^*(\alpha)\Big)< \beta < \gamma < \min\Big(2,\frac\pi\alpha\Big)  ,
\]
$\gamma>0$ and $\beta\not=0$. Furthermore, it is assumed that the integral of $g$ over $K$ is zero. Then there exist solutions
$(u,p)\in E_\gamma^2(K)\times V_\gamma^1(K)$ and $(u',p')\in E_\beta^2(K)\times V_\beta^1(K)$ of the problem (\ref{2bt2}).
By \cite[Lemma 4.5 and Lemma 4.7]{r-18}, we have $u=u'$ and, consequently, $\nabla p = \nabla p'$. If $\beta>0$ then even $p=p'$.
Thus, the solution $(u,p)$ satisfies the estimate
\[
\| u\|_{V_\beta^2(K)} + |s|\, \| u\|_{V_\beta^0(K)} + \| \nabla p\|_{V_\beta^0(K)} \le c\, \Big( \| f\|_{V_\beta^0(K)} + \| g\|_{V_\beta^1(K)}
 + |s|\, \| g\|_{(W_{-\beta}^1(K))^*} \Big)
\]
with a constant $c$ independent of $f$, $g$ and $s$.}
\end{Rem}

\section{The parameter-dependent problem in a polygon}

We consider the problem (\ref{par3}), (\ref{par4}) in the polygonal domain $\Omega$ which was described in the introduction.
First, we prove the existence of solutions $(u,p) \in V_\beta^2(\Omega) \times W_\beta^1(\Omega)$.
After this, we estimate this solution in relation to the data $f,g$ and the parameter $s$.

\subsection{Existence of solutions}

Let $W^l(\Omega)$ be the Sobolev space with the norm
\[
\| u\|_{W^l(\Omega)} = \Big( \int_\Omega \sum_{|\alpha|\le l} |\partial_x^\alpha u(x)|^2\, dx\Big)^{1/2} .
\]
The closure of the set $C_0^\infty(\Omega)$ in this space is denoted by $\stackrel{\circ}W\!{}^l(\Omega)$, and
$W^{-l}(\Omega)$ is defined as the dual space of $\stackrel{\circ}W\!{}^l(\Omega)$. By $(\cdot,\cdot)_\Omega$ we denote
the $L_2(\Omega)$-scalar product and its extension to the product $W^{-l}(\Omega)\times\stackrel{\circ}W\!{}^l(\Omega)$.

Furthermore, let $\stackrel{\circ}{L}{}_2(\Omega)$ be the set of all $g\in L_2(\Omega)$ satisfying the condition
\begin{equation} \label{1cl1}
\int_\Omega g(x)\, dx =0.
\end{equation}
We consider the bilinear form
\[
b_s(u,v) = \int_\Omega \Big( su\cdot v + 2\sum_{i,j=1}^3 \varepsilon_{i,j}(u)\, \varepsilon_{i,j}(v)\Big)\, dx,
\]
where $\varepsilon(u)$ denotes the strain tensor with the elements
\[
\varepsilon_{i,j}({\mathfrak u})=\frac 12\,\Big( \frac{\partial u_i}{\partial x_j}+\frac{\partial u_j}{\partial x_i}\Big), \quad
i,j=1,\ldots,2,
\]
Obviously, this bilinear form is continuous on $W^1(\Omega)\times W^1(\Omega)$. Furthermore, it follows from Korn's inequality  that
\begin{equation} \label{2cl1}
\big| b_s(u,\overline{u})\big| \ge c\, \|u\|^2_{W^1(\Omega)}
\end{equation}
for all $u\in W^1(\Omega)$ and for all $s$, $\mbox{Re}\, s\ge 0$, $s\not=0$, where $c$ depends on $|s|$ but not on $u$.
The following lemma can be found e.g. in \cite[Chapter 1, Corollary 2.4]{Girault} or \cite[Lemma 11.1.1]{mr-10}

\begin{Le} \label{cl1}
Let $g\in \stackrel{\circ}{L}\!{}_2(\Omega)$. Then there exists a vector function $u\in \stackrel{\circ}W\!{}^1(\Omega)$
such that $\nabla\cdot u=g$ in $\Omega$ and
\[
\| u\|_{W^1(\Omega)} \le c\, \| g\|_{L_2(\Omega)}
\]
with a constant $c$ independent of $g$.
\end{Le}

Applying the inequality (\ref{2cl1}) and Lemma \ref{cl1}, one can prove the following lemma analogously to \cite[Theorem 11.1.2]{mr-10}.

\begin{Le} \label{cl2}
Suppose that $f\in W^{-1}(\Omega)$, $g\in \stackrel{\circ}{L}\!{}_2(\Omega)$,
$\mbox{\em Re}\, s\ge 0$ and $s\not=0$. Then there exists a unique solution $(u,p) \in \stackrel{\circ}W\!{}^1(\Omega)\times \stackrel{\circ}{L}\!{}_2(\Omega)$
of the problem
\begin{equation} \label{1cl2}
b_s(u,\bar{v}) - \int_\Omega p\, \nabla\cdot\bar{v}\, dx = (f,v)_\Omega \ \mbox{ for all }v\in \stackrel{\circ}W\!{}^1(\Omega),\quad
-\nabla\cdot u =g\mbox{ in }\Omega.
\end{equation}
satisfying the estimate
\begin{equation} \label{2cl2}
\| u\|_{W^1(\Omega)} + \| p\|_{L_2(\Omega)} \le c\, \Big( \| f\|_{W^{-1}(\Omega)} + \| g\|_{L_2(\Omega)}\Big),
\end{equation}
where $c$ depends on $|s|$ but not on $f$ and $g$.
\end{Le}

The solution of the problem (\ref{1cl2}) is called a weak solution of the problem (\ref{par3}), (\ref{par4}).
Let $l$ be a positive integer. Then we define $W_\beta^l(\Omega)$ as the weighted Sobolev space with the norm
\[
\| u\|_{W_\beta^l(\Omega)} = \Big( \| \nabla u\|^2_{V_\beta^{l-1}(\Omega)} + \| u \|^2_{L_2(\Omega)}\Big)^{1/2}.
\]
It follows from Hardy's inequality that $W_\beta^l(\Omega)=V_\beta^l(\Omega)$ if $\beta_j>l-1$ for all $j$.
Using the last lemma together with well-known regularity results for elliptic problems (and, in particular, for the stationary Stokes system),
we are able to prove the following lemma.

\begin{Le} \label{cl3}
Suppose that $\mbox{\em Re}\, s \ge 0$, $s\not=0$, $f\in V_\beta^0(\Omega)$ and $g\in V_\beta^1(\Omega)$,
where the components $\beta_j$ of $\beta$ satisfy the inequalities
\[
1-\mbox{\em Re}\, \lambda^*(\alpha_j)< \beta_j <1, \ \ \beta_j\not=0 \ \mbox{ for all }j.
\]
Furthermore, we assume that $g$ satisfies the condition {\em (\ref{1cl1})}.
Then there exists a unique solution $(u,p)\in V_\beta^2(\Omega)\times (W_\beta^1(\Omega)\cap\stackrel{\circ}{L}\!{}_2(\Omega))$
of the problem {\em (\ref{par3}), (\ref{par4})} satisfying the estimate
\begin{equation} \label{1cl3}
\| u \|_{V_\beta^2(\Omega)} + \|  p\|_{W_\beta^1(\Omega)} \le c\, \Big( \| f\|_{V_\beta^0(\Omega)} + \| g\|_{V_\beta^1(\Omega)}\Big),
\end{equation}
where $c$ depends on $|s|$ but not on $f$ and $g$.
\end{Le}

Proof. Since $\beta_j<1$ for all $j$, it follows from H\"older's and Hardy's inequalities that
\[
\Big| \int_\Omega f\cdot v\, dx\Big| \le \| f\|_{V_\beta^0(\Omega)} \ \| v\|_{V_{-\beta}^0(\Omega)}
  \le c\, \| f\|_{V_\beta^0(\Omega)} \ \| v\|_{W^1(\Omega)}
\]
Hence, the functional
\[
\stackrel{\circ}{W}\!{}^1(\Omega) \ni v \to (f,\bar{v})_\Omega = \int_\Omega f\cdot v\, dx
\]
is continuous. Furthermore, $g\in V_\beta^1(\Omega)\subset L_2(\Omega)$. By Lemma \ref{cl2}, there exists a unique solution
$(u,p)\in \stackrel{\circ}{W}\!{}^1(\Omega) \times \stackrel{\circ}{L}\!{}_2(\Omega)$ of the problem (\ref{1cl2}) satisfying the estimate (\ref{2cl2}).
Let $\zeta_j$ be two times continuously differentiable functions with support in a sufficiently small neighborhood of $P_j$ equal to one
near $P_j$. The vector function $\zeta_j(u,p)$ satisfies the equations
\[
-\Delta(\zeta_ju) + \nabla(\zeta_j p) = F_j
  = \zeta_j(f-\nabla g-su) -u\, \Delta \zeta_j - 2\sum_{i=1}^2 \partial_{x_i}\zeta_j \, \partial_{x_i} u + p\, \nabla\zeta_j
\]
and $-\nabla \cdot (\zeta_j u) = G_j = \zeta_j g -u\cdot \nabla\zeta_j$. It follows from Hardy's inequality that $\zeta_j u \in V_{\beta_j}^0(K_j)$
for $u\in \stackrel{\circ}{W}\!{}^1(\Omega)$ and $\beta_j>-1$. Consequently, $F_j\in V_{\beta_j}^0(K_j)$ and $G_j \in V_{\beta_j}^1(K_j)$ if $\beta_j >-1$.
Using well-known regularity results for the stationary Stokes system, we conclude that $\zeta_j (u,p) \in V_\beta^2(\Omega)\times V_\beta^1(\Omega)$
if $\beta_j>0$ and
\[
\| \zeta_j u\|_{V_\beta^2(\Omega)} + \|\zeta_j p\|_{V_\beta^1(\Omega)} \le c\, \Big( \| F_j\|_{V_{\beta_j}^0(K_j)} + \| G_j \|_{V_{\beta_j}^1(K_j)}\Big)
  \le c' \, \Big( \| f\|_{V_\beta^0(\Omega)} + \| g\|_{V_\beta^1(\Omega)}\Big).
\]
If $\max(-1,1-\mbox{Re}\, \lambda^*(\alpha_j))<\beta_j<0$, then the vector function $\zeta_j(u,p)$ admits the asymptotics $\zeta_j(u,p)=(0,c_j) + (v,q)$,
where $(v,q)\in V_{\beta_j}^2(K_j)\times V_{\beta_j}^1(K_j)$.
This implies $\zeta_j u\in V_\beta^2(\Omega)$ and $\nabla (\zeta_j p)\in  V_\beta^0(\Omega)$. Furthermore, the estimate (\ref{1cl3}) holds.

Suppose that $\max(-3,1-\mbox{Re}\, \lambda^*(\alpha_j))< \beta_j \le -1$. Then $F_j \in V_{\beta'_j}^0(K_j)$ and $G_j \in V_{\beta'_j}^0(K_j)$
for arbitrary $\beta'_j >-1$. Consequently, $\zeta_j u\in V_{\beta'_j}^2(\Omega)$ and $\nabla (\zeta_j p)\in  V_{\beta'_j}^0(\Omega)$ for arbitrary
$\beta'_j$, $-1<\beta'_j>0$.
But then $F_j\in V_{\beta_j}^0(K_j)$ and $G_j \in V_{\beta_j}^1(K_j)$ and, consequently, $\zeta_j u\in V_\beta^2(\Omega)$ and $\nabla (\zeta_j p)\in  V_\beta^0(\Omega)$.

Repeating this argument, we obtain this result for arbitrary $\beta_j>1-\mbox{Re}\, \lambda^*(\alpha_j)$. This proves the existence of a solution
$(u,p)\in V_\beta^2(\Omega)\times (W_\beta^1(\Omega)\cap\stackrel{\circ}{L}\!{}_2(\Omega))$. The uniqueness of this solution follows from the
imbedding $V_\beta^2(\Omega) \subset W^1(\Omega)$  (since $\beta_j<1$ for all $j$) and Lemma \ref{cl2}.  $\Box$

\subsection{An estimate for \boldmath $p$\unboldmath}

Clearly the constant $c$ in the estimate  of Lemma \ref{cl3} depends on $s$. Our goal is to obtain a more precise estimate.
We prove the existence  of solutions of the Neumann problem for the Poisson equation in the function space $W_\beta^2(\Omega)$ which was introduced in the foregoing subsection.

\begin{Le} \label{cl4}
Suppose that $\phi \in V_{1-\gamma}^0(\Omega)$, where $0<\gamma_j < \pi/\alpha_j$ for all $j$, and that the integral of $\phi$ over $\Omega$ is zero.
Then there exists a solution $q\in W_{1-\gamma}^2(\Omega)$ of the problem
\begin{equation} \label{1cl4}
-\Delta q = \phi \ \mbox{ in }\Omega, \quad \frac{\partial q}{\partial n}=0\ \mbox{on }\Gamma
\end{equation}
which satisfies the estimate
\begin{equation} \label{2cl4}
\| q\|_{L_2(\Omega)} + \| \nabla q\|_{V_{1-\gamma}^1(\Omega)} \le c \, \| \phi \|_{V_{1-\gamma}^0(\Omega)}
\end{equation}
with a constant $c$ independent of $\phi$.
\end{Le}

Proof. Since $\gamma_j>0$ for all $j$, H\"older's inequality implies the imbedding $V_{1-\gamma}^0(\Omega) \subset L_1(\Omega)$.
Furthermore, it follows from Hardy's inequality that
\[
\int_\Omega \prod r_j^{2\gamma_j-2} |v|^2\, dx \le c \int_\Omega \big( \prod r_j^{2\gamma_j}\, |\nabla v|^2 + |v|^2\big)\, dx \le c'\, \| v\|^2_{W^1(\Omega)}
\]
for all $v \in W^1(\Omega)$, i.e., $W^1(\Omega) \subset V_{\gamma-1}^0(\Omega)$ and, consequently,
$V_{1-\gamma}^0(\Omega) \subset (W^1(\Omega))^*$.
As is known, the $W^1$-norm is equivalent to the norm
\[
\| q\| = \| \nabla q\|_{L_2(\Omega)}
\]
on the subspace $W^1(\Omega) \cap \stackrel{\circ}{L}\!{}_2(\Omega)$.
Hence there exists a unique variational solution $q\in W^1(\Omega) \cap \stackrel{\circ}{L}\!{}_2(\Omega)$ of the problem
(\ref{1cl4}), i.e.,
\[
\int_\Omega \nabla q\cdot \nabla v\, dx = \int_\Omega \phi\, v\, dx \ \mbox{ for all }v\in W^1(\Omega).
\]
This solution satisfies the estimate
\[
\| q\|_{W^1(\Omega)} \le c'\, \| \phi\|_{(W^1(\Omega))^*} \le c\, \| \phi\|_{V_{1-\gamma}^0(\Omega)}\, .
\]
For every $j$, let $\zeta_j$ be a smooth (of class $C^2$) cut-off function with support in the neighborhood ${\cal U}_j$ of $P_j$
which is equal to one near $P_j$ and satisfies the condition $\partial \zeta_j/\partial n =0$ on $\Gamma$.
Obviously, $\zeta_j q \in V_\varepsilon^1(K_j)$ with arbitrary positive $\varepsilon$ and
\[
-\Delta(\zeta_j q) = \zeta_j \phi -2\nabla\zeta_j\cdot \nabla q - q\, \Delta\zeta_j \in V_{1-\gamma_j}^0(K_j).
\]
Hence, $\zeta_j q$ admits the representation
\[
\zeta_j q = c_j+ d_j\log r_j + w_j,
\]
where $w_j\in V_{1-\gamma_j}^2(K_j)$. Since $\zeta_j(q-w_j) \in W^1(K_j)$, it follows that $d_j=0$ and, consequently,  $\nabla(\zeta_j q)
= \nabla w_j\in V_{1-\gamma_j}^1(K_j)$,
\[
\| \nabla(\zeta_j q)\|_{V_{1-\gamma_j}^1(K_j)} \le \| w\|_{V_{1-\gamma_j}^2(K_j)} \le c\, \| \phi\|_{V_{1-\gamma}^0(\Omega)}\, .
\]
Since obviously, $q\in W^2(\Omega_\varepsilon)$ for every subdomain $\Omega_\varepsilon=\{ x\in \Omega:\ \mbox{dist}(x,{\cal P})>\varepsilon\}$,
we conclude that $\nabla q\in V_{1-\gamma}^1(\Omega)$. Furthermore, the estimate (\ref{2cl4}) holds. $\Box$ \\

In the following lemma, we consider solutions $(u,p) \in V_\gamma^2(\Omega)\times V_\gamma^1(\Omega)$ of the problem (\ref{par3}), (\ref{par4}).
Note that the constant function $p=c$ is an element of the space $V_\gamma^1(\Omega)$ if $\gamma_j>0$ for all $j$.
As in Lemma \ref{cl3}, we assume that the integral of $p$ over $\Omega$ is zero. This integral exists for every $p\in V_\gamma^1(\Omega)$ if
$\gamma_j<2$ for all $j$.

\begin{Le} \label{cl6}
Suppose that $\mbox{\em Re}\, s\ge 0$, $s\not=0$, and that $(u,p) \in V_\gamma^2(\Omega)\times V_\gamma^1(\Omega)$
is a solution of the problem {\em (\ref{par3}), (\ref{par4})}. We assume that
\[
0< \gamma_j < \min\big(2,\frac\pi{\alpha_j}\big)
\]
for all $j$ and that the integral of $p$ over $\Omega$ is zero. Then
\begin{eqnarray} \label{1cl6}
\| p\|_{V_{\gamma-1}^0(\Omega)}  \le  c\, \Big( \| f\|_{V_\gamma^0(\Omega)} + \| g\|_{V_\gamma^1(\Omega)}
  + |s|\, \| g\|_{(W_{-\gamma}^1(\Omega))^*}
 +  \| Du\|_{V_{\gamma-1/2}^0(\Gamma)}\Big),
\end{eqnarray}
where $Du$ is the matrix with the elements $\partial_{x_i} u_j$ and
\[
\| u\|_{V_{\gamma-1/2}^0(\Gamma)} = \Big\| \prod_j r_j^{\gamma_j-1/2} u\Big\|_{L_2(\Gamma)}\, .
\]
The constant $c$ in {\em (\ref{1cl6})} is independent of $u,p$ and $s$.
\end{Le}

Proof.
By (\ref{par3}), (\ref{par4}), we have
\[
\int_\Omega \nabla p\cdot \nabla \bar{q}\, dx = (\Phi,q)_\Omega \ \mbox{for all } q\in W_{1-\gamma}^2(\Omega),
\]
where
\[
(\Phi,q)_\Omega = \int_\Omega (f+\Delta u+\nabla\nabla\cdot u)\cdot\nabla\bar{q}\, dx  -s\, (g,q)_\Omega \, .
\]
Integration by parts yields
\begin{equation} \label{2cl6}
-\int_\Omega p\, \Delta \bar{q}\, dx + \int_\Gamma p\, \frac{\partial \bar{q}}{\partial n}\, d\sigma = (\Phi,q)_\Omega.
\end{equation}
Let $\phi \in V_{1-\gamma}^0(\Omega)$. By H\"older's inequality,
\[
\Big| \int_\Omega \phi\, dx \Big|\le c\, \| \phi\|_{V_{1-\gamma}^0(\Omega)}\, .
\]
We define $\displaystyle c_0=\frac 1{|\Omega|} \int_\Omega \phi\, dx.$ By Lemma \ref{cl4},
there exists a solution $q\in W_{1-\gamma}^2(\Omega)$ of the problem
\[
-\Delta q = \phi-c_0 \ \mbox{ in }\Omega, \quad \frac{\partial q}{\partial n} =0 \ \mbox{ on }\Gamma
\]
satisfying the estimate
\begin{equation} \label{4cl6}
\| \nabla q\|_{V_{-\gamma}^0(\Omega)} \le \| q\|_{W_{1-\gamma}^2(\Omega)} \le c\, \| \phi-c_0\|_{V_{1-\gamma}^0(\Omega)} \le c'\, \| \phi\|_{V_{1-\gamma}^0(\Omega)}\, .
\end{equation}
Using (\ref{2cl6}), we obtain
\begin{equation} \label{3cl6}
 \int_\Omega p\, \bar{\phi}\, dx =  \int_\Omega p\, \overline{\phi-c_0}\, dx = -\int_\Omega p\, \Delta \bar{q}\, dx = (\Phi,q)_\Omega  \, .
\end{equation}
We set $\phi = r^{2\gamma-2}p$ and obtain
\[
\| p\|^2_{V_{\gamma-1}^0(\Omega)} = (\Phi,q)_\Omega  \, .
\]
There is the decomposition $\Phi=\Phi_1+\Phi_2$, where
\[
(\Phi_1,q)_\Omega = \int_\Omega \Big( (f-2\nabla g)\cdot\nabla\bar{q} -sg\,\bar{q}\Big)\, dx, \quad
(\Phi_2,q)_\Omega = \int_\Omega ( \Delta u-\nabla\nabla\cdot u)\, \nabla\bar{q}\, dx
\]
Obviously,
\[
\Big| \int_\Omega (f-2\nabla g)\, \nabla\bar{q}\, dx\Big| \le c\, \Big( \| f\|_{V_\gamma^0(\Omega)} + \| g\|_{V_\gamma^1(\Omega)}\Big)
  \, \|\nabla q\|_{V_{-\gamma}^0(\Omega)}\, .
\]
and
\[
\Big| \int_\Omega g\, \bar{q}\, dx\Big| \le  \| g\|_{(W_{-\gamma}^1(\Omega))^*} \, \| q\|_{W_{-\gamma}^1(\Omega)} \le
  \| g\|_{(W_{-\gamma}^1(\Omega))^*} \, \| q\|_{W_{1-\gamma}^2(\Omega)}\, .
\]
Thus,
\[
\big| (\Phi_1,q)_\Omega\big| \le c\, \Big( \| f\|_{V_\gamma^0(\Omega)} + \| g\|_{V_\gamma^1(\Omega)}
  + |s|\, \| g\|_{(W_{-\gamma}^1(\Omega))^*}\Big)\, \|p\|_{V_{\gamma-1}^0(\Omega)}\, .
\]
Furthermore,
\begin{eqnarray*}
&& \big| (\Phi_2,q)_\Omega\big| = \Big| \int_\Gamma \sum_{i,j=1}^2 \frac{\partial u_i}{\partial x_j}
\Big( n_j\, \frac{\partial \bar{q}}{\partial x_i} -n_i\, \frac{\partial\bar{q}}{\partial x_j}\Big)\, d\sigma \Big| \\
&& \le  c\, \| Du\|_{V_{\gamma-1/2}^0(\Gamma)}\, \| \nabla q\|_{V_{-\gamma+1/2}^0(\Gamma)}
 \le  c'\, \| Du\|_{V_{\gamma-1/2}^0(\Gamma)} \, \| \nabla q\|_{V_{1-\gamma}^1(\Omega)}\, ,
\end{eqnarray*}
Consequently, (\ref{3cl6}) together with (\ref{2cl4}) yields (\ref{1cl6}). $\Box$

\subsection{An estimate for the solution}

Now, we prove the main result of this section.

\begin{Th} \label{ct1}
Suppose that $\mbox{\em Re}\, s \ge 0$, $s\not=0$, $f\in V_\beta^0(\Omega)$ and $g\in V_\beta^1(\Omega)$, where the components
$\beta_j$ of $\beta$ satisfy the inequalities {\em (\ref{3ct1})} for all $j$. Furthermore, we assume that $g$ satisfies the condition {\em (\ref{1cl1})}.
Then there exists a unique solution $(u,p) \in V_\beta^2(\Omega) \times (W_\beta^1(\Omega)\cap \stackrel{\circ}{L}\!\!{}_2(\Omega))$
of the problem {\em (\ref{par3}), (\ref{par4})}. Moreover, there exists a positive number $\gamma_0$ such that
the solution $(u,p)$  satisfies the estimate
\begin{equation} \label{1ct1}
\| u\|^2_{V_\beta^2(\Omega)} + |s|^2\, \| u \|^2_{V_\beta^0(\Omega)} + \| p\|^2_{W_\beta^1(\Omega)} \le  c\, \Big( \| f\|^2_{V_\beta^0(\Omega)}
+ \| g\|^2_{V_\beta^1(\Omega)}  + \ |s|^2\, \| g\|^2_{(W_{-\beta}^1(\Omega))^*} \Big)
\end{equation}
for $\mbox{\em Re}\, s \ge 0$, $|s|>\gamma_0$. Here, the constant $c$ is independent of $f,g$ and $s$.
\end{Th}

Proof.
The existence and uniqueness of a solution $(u,p)\in V_\beta^2(\Omega) \times (W_\beta^1(\Omega)\cap \stackrel{\circ}{L}\!\!{}_2(\Omega))$
follows from Lemma \ref{cl3}. We prove the estimate (\ref{1ct1}).

Let $\zeta_1,\ldots,\zeta_n$ be the same cut-off functions as in the proof of Lemma \ref{cl3}.
Furthermore, let $\zeta_{n+1}$ be such that $\zeta_1+\cdots+\zeta_{n+1}=1$ in $\Omega$.
Then the vector function $\zeta_j(u,p)$ satisfies the equations
\[
(s-\Delta)(\zeta_j u)-\nabla\nabla\cdot(\zeta_j u)+\nabla(\zeta_j p)= F_j,\quad -\nabla\cdot (\zeta_j u)=G_j\ \mbox{ in }\Omega,
\]
where $F_j = \zeta_j( f+\Delta u+\nabla\nabla\cdot u) -\Delta(\zeta_j u)  - \nabla\nabla\cdot(\zeta_j u)+p\nabla\zeta_j$
and $G_j=\zeta_j g-u\cdot\nabla\zeta_j$.
Furthermore, $\zeta_j u =0$ on $\partial\Omega$.
By Theorem \ref{bt2} and Remark \ref{br2}, there is the estimate
\begin{eqnarray} \label{2ct1} \nonumber
&& \| \zeta_j u\|_{V_\beta^2(\Omega)} + |s|\, \| \zeta_j u\|_{V_\beta^0(\Omega)} + \| \nabla(\zeta_j p)\|_{V_\beta^0(\Omega)} \\
&& \le c\, \Big( \| F_j\|_{V_\beta^0(\Omega)} + \| G_j \|_{V_\beta^1(\Omega)} + |s|\, \| G_j\|_{(W_{-\beta}^1(\Omega))^*}\Big)
\end{eqnarray}
for $j=1,\ldots,n$. Since all functions in (\ref{2ct1}) have supports in a neighborhood of $P_j$ for $j=1,\ldots,n$, their norms in
$V_\beta^l(\Omega)$ are equivalent to the $V_{\beta_j}^l(K_j)$-norms. The function $\zeta_{n+1}$ is zero in a neighborhood of any corner point
$P_1,\ldots,P_n$. Using existence and uniqueness results for the Dirichlet problem for the nonstationary
Stokes system in domains with smooth boundaries (see \cite[Theorem 3.1]{sol-03}), one can prove the estimate (\ref{2ct1}) for $j=n+1$.
Here, we refer to \cite[Lemma 2.2]{kr-16} and \cite[Lemma 2.6]{kr-20} (in \cite{kr-16}, the norm of $g$ in $W^{-1}$ must be replaced by the norm in $(W^1)^*$).
One can easily verify the estimate
\begin{eqnarray} \label{3ct1a} \nonumber
&& \| F_j\|_{V_\beta^0(\Omega)} + \| G_j \|_{V_\beta^1(\Omega)} + |s|\, \| G_j\|_{(W_{-\beta}^1(\Omega))^*} \le
  \| \zeta_j f\|_{V_\beta^0(\Omega)} + \| \zeta_j g \|_{V_\beta^1(\Omega)} \\
&& \quad + |s|\, \| \zeta_j g\|_{(W_{-\beta}^1(\Omega))^*}
 + c\, \Big(\| u\|_{W^1(\Omega_j)} + \| p\|_{L_2(\Omega_j)} + |s|\, \| u\cdot \nabla\zeta_j\|_{(W_{-\beta}^1(\Omega))^*}\Big),
\end{eqnarray}
where $\Omega_j = \Omega \cap \mbox{supp}\, \nabla\zeta_j$.
By Ehrling's lemma, there is the inequality
\[
\| u\|_{W^1(\Omega_j)} \le \varepsilon\, \| u\|_{W^2(\Omega_j)} + c(\varepsilon)\, \| u\|_{L_2(\Omega_j)}\le
c_1\, \varepsilon\, \| u\|_{V_\beta^2(\Omega)} + c_1\, c(\varepsilon)\, \| u\|_{V_\beta^0(\Omega)}
\]
with an arbitrarily small positive $\varepsilon$. We estimate the norm of $u\cdot\nabla\zeta_j$ in $(W_{-\beta}^1(\Omega))^*$. Let $q\in W_{-\beta}^1(\Omega)$. Then
\begin{eqnarray*}
&& \int_{\Omega} sq u \cdot \nabla\zeta_j\, dx = \int_{\Omega} (f-\nabla g+\Delta u -\nabla p)\cdot q\nabla\zeta_j\, dx
 =  \int_{\Omega}  (f-\nabla g)\cdot q\nabla\zeta_j\, dx \\
&& \qquad - \int_\Omega\Big( \sum_{i=1}^2\nabla u_i\cdot \nabla(q\partial_{x_i}\zeta_j)+ p\,\nabla\cdot (q\nabla\zeta_j)\Big)\, dx
   + \int_\Gamma \Big( \frac{\partial u}{\partial n}\cdot q\nabla\zeta_j - pq\nabla\zeta_j\cdot n\Big)\, d\sigma
\end{eqnarray*}
and, consequently,
\begin{eqnarray*}
&& \| su\cdot\nabla\zeta_j \|_{(W_{-\beta}^1(\Omega))^*} \\
&& \le c\, \Big( \| f-\nabla g\|_{V_\beta^0(\Omega)}+ \| u\|_{W^1(\Omega_j)}
  + \| p\|_{L_2(\Omega_j)} + \Big\| \frac{\partial u}{\partial n}\cdot \nabla\zeta_j\Big\|_{L_2(\Gamma)} 
  + \Big\| p\, \frac{\partial \zeta_j}{\partial n}\Big\|_{L_2(\Gamma)}\Big).
\end{eqnarray*}
Here,
\begin{eqnarray*}
\Big\| \frac{\partial u}{\partial n}\cdot \nabla\zeta_j\Big\|_{L_2(\Gamma)} + \Big\| p\, \frac{\partial \zeta_j}{\partial n}\Big\|_{L_2(\Gamma)}
\le c\, \Big( \| u\|_{W^{1+\sigma}(\Omega_j)} + \| p\|_{W^\sigma(\Omega_j)}\Big)
\end{eqnarray*}
with an arbitrary $\sigma \in (\frac 12,1)$. Using Ehrling's lemma, we obtain
\begin{eqnarray*}
&& \Big\| \frac{\partial u}{\partial n}\cdot \nabla\zeta_j\Big\|_{L_2(\Gamma)} + \Big\| p\, \frac{\partial \zeta_j}{\partial n}\Big\|_{L_2(\Gamma)}\\
&& \le  \varepsilon \,\Big(  \| u\|_{V_\beta^2(\Omega)} + \| p\|_{W_\beta^1(\Omega)}\Big)
 +  c(\varepsilon)\, \Big(  \| u\|_{V_\beta^0(\Omega)} + \| p\|_{L_2(\Omega)}\Big).
\end{eqnarray*}
It remains to estimate the norm of $p$ in $L_2(\Omega)$. Let $\gamma_j$ be real numbers $\max(0,\beta_j)<\gamma_j<\min(1,\frac{\pi}{\alpha_j})$
for $j=1,\ldots,n$.
By Lemma \ref{cl6}, $p$ satisfies the estimate
\begin{eqnarray*}
\| p\|_{L_2(\Omega)} & \le & c \, \| p\|_{V_{\gamma-1}^0(\Omega)} \\
& \le & c'\, \Big( \| f\|_{V_\gamma^0(\Omega)} + \| g\|_{V_\gamma^1(\Omega)}
   + |s|\, \| g\|_{(W_{-\gamma}^1(\Omega))^*} + \| Du\|_{V_{\gamma-1/2}^0(\Gamma)}\Big).
\end{eqnarray*}
Let $\varepsilon'$ be an arbitrarily small positive number. Since $V_\beta^{1/2}(\Gamma) \subset V_{\beta-1/2}^0(\Gamma)$
and $\beta_j<\gamma_j$ for all $j$, one can choose a subset $\Gamma'$ of $\Gamma$ with positive distance to ${\cal P}$
such that
\begin{eqnarray*}
\| Du\|_{V_{\gamma-1/2}^0(\Gamma)} & \le & \varepsilon' \, \| Du\|_{V_\beta^{1/2}(\Gamma)} + c(\varepsilon')\, \| Du\|_{L_2(\Gamma')} \\
& \le & \varepsilon' \, \| Du\|_{V_\beta^1(\Omega)} + c(\varepsilon')\, \Big( \varepsilon''\, \| Du\|_{V_\beta^1(\Omega)}
  + c_1(\varepsilon'')\, \| u\|_{V_\beta^0(\Omega)}\Big),
\end{eqnarray*}
where $\varepsilon''$ can be chosen arbitrarily small. Hence,
\[
\| p\|_{L_2(\Omega)}  \le  c\, \Big( \| f\|_{V_\gamma^0(\Omega)} + \| g\|_{V_\gamma^1(\Omega)}
   + |s|\, \| g\|_{(W_{-\gamma}^1(\Omega))^*}\Big) + \varepsilon' \, \| u\|_{V_\beta^2(\Omega)} + C(\varepsilon')\, \| u\|_{V_\beta^0(\Omega)}\, ,
\]
where $\varepsilon'$ can be chosen arbitrarily small.
Thus, the estimates (\ref{2ct1}), (\ref{3ct1a}) together with the above estimates for the $W^1$-norm of $u$, the $L_2$-norm of $p$ on $\Omega_j$ and
the norm of $u\cdot \nabla\zeta_j$ lead to the estimate
\begin{eqnarray*}
&& \| u\|_{V_\beta^2(\Omega)} + |s|\, \| u \|_{V_\beta^0(\Omega)} + \| p\|_{W_\beta^1(\Omega)} \le  c\, \Big( \| f\|_{V_\beta^0(\Omega)}
+ \| g\|_{V_\beta^1(\Omega)}  + \ |s|\, \| g\|_{(W_{-\beta}^1(\Omega))^*} \Big) \\
&& \quad  + \frac 12\, \Big( \| u\|_{V_\beta^2(\Omega)}
+ \| p\|_{W_\beta^1(\Omega)}\Big) + C\, \| u\|_{V_\beta^0(\Omega)}
\end{eqnarray*}
with constants $c$ and $C$ independent of $u,p$ and $s$. For $|s|>2C$, the inequality (\ref{1ct1}) holds. $\Box$

\section{The time-dependent problem on \boldmath $\Omega$}

Now, we consider the problem (\ref{stokes1}), (\ref{stokes2}).

\subsection{Weighted Sobolev spaces in $\Omega\times (0,T)$}

Let $0<T\le \infty$ and $l$ be a nonnegative integer. Then $L_2(0,T;V_\beta^l(\Omega))$ is defined as the space of all functions (vector functions)
on $\Omega\times (0,T))$ with finite norm
\[
\| {\mathfrak u}\|_{L_2(0,T;V_\beta^l(\Omega))} = \Big( \int_0^T  \| {\mathfrak u}(\cdot,t) \|^2_{V_\beta^l(\Omega)}\, dt \Big)^{1/2}.
\]
Furthermore, let
$\stackrel{\circ}{W}\!{}_\beta^{2,1}(\Omega\times (0,T))$ be the space of all functions (vector functions)
${\mathfrak u}(x,t)$ on $\Omega\times (0,T)$ with finite norm
\[
\| {\mathfrak u}\|_{W_\beta^{2,1}(\Omega\times (0,T))} = \Big( \int_0^T  \big( \| {\mathfrak u}(\cdot,t) \|^2_{V_\beta^2(\Omega)}
  + \| \partial_t {\mathfrak u}(\cdot,t) \|^2_{V_\beta^0(\Omega)}\big)\, dt \Big)^{1/2}
\]
which are zero for $t=0$.
Analogously, $\stackrel{\circ}{W}\!{}_\beta^{1,1}(\Omega\times (0,T))$ is defined as the space of all functions (vector functions)
on $\Omega\times (0,T)$ with finite norm
\[
\| {\mathfrak u}\|_{W_\beta^{1,1}(\Omega\times (0,T))} = \Big( \int_0^T  \big( \| {\mathfrak u}(\cdot,t) \|^2_{V_\beta^1(\Omega)}
  + \| \partial_t {\mathfrak u}(\cdot,t) \|^2_{(W_{-\beta}^1(\Omega))^*}\big)\, dt \Big)^{1/2}
\]
which vanish for $t=0$.

Moreover, we define $\stackrel{\circ}{W}\!{}_\beta^{l,1}(\Omega\times {\Bbb R}_+,e^{-\gamma t})$ for $l=1,2$ and
$L_{2,-\gamma}(0,\infty;V_\beta^l(\Omega))$ as the sets of functions ${\mathfrak u}={\mathfrak u}(x,t)$ such that
$e^{-\gamma t} {\mathfrak u}\in\stackrel{\circ}{W}\!{}_\beta^{l,1}(\Omega\times (0,\infty))$ and $e^{-\gamma t} {\mathfrak u}\in L_2(0,\infty;V_\beta^l(\Omega))$,
respectively. The space $\stackrel{\circ}{W}\!{}_\beta^{l,1}(\Omega\times {\Bbb R}_+,e^{-\gamma t})$ is provided with the norm
\[
\| {\mathfrak u} \|_{W_\beta^{l,1}(\Omega\times {\Bbb R}_+,e^{-\gamma t})}
  = \| e^{-\gamma t}{\mathfrak u} \|_{W_\beta^{1,1}(\Omega\times (0,\infty))}\, .
\]
Analogously, the norm in $L_{2,-\gamma}(0,\infty;V_\beta^l(\Omega))$ is defined.

We consider the Laplace transforms. Let $H_\beta(\Omega,\gamma)$ be the space of holomorphic functions
$u(x,s)$ for $\mbox{Re}\, s>\gamma$ with values in $V_\beta^0(\Omega)$ for which the norm
\[
\| u\|_{H_\beta(\Omega,\gamma)} = \sup_{\sigma>\gamma} \Big( \frac 1i \int_{\mbox{\scriptsize Re}\, s=\sigma} \| u(\cdot,s)\|^2_{V_\beta^0(\Omega)}\, ds  \Big)^{1/2}
\]
is finite. The spaces $H_\beta^l(\Omega,\gamma)$, $l=1,2$, are the sets of holomorphic functions
$u(x,s)$ for $\mbox{Re}\, s>\gamma$ with values in $V_\beta^l(\Omega)$ for which the norms
\[
\| u\|_{H_\beta^1(\Omega,\gamma)} = \sup_{\sigma>\gamma} \Big( \frac 1i \int_{\mbox{\scriptsize Re}\, s=\sigma}  \big( \| u(\cdot,s)\|^2_{V_\beta^1(\Omega)}
  + |s|^2  \| u(\cdot,s)\|^2_{(W_{-\beta}^1(\Omega))^*}\big)\, ds  \Big)^{1/2}
\]
and
\[
\| u\|_{H_\beta^2(\Omega,\gamma)} = \sup_{\sigma>\gamma} \Big( \frac 1i \int_{\mbox{\scriptsize Re}\, s=\sigma}  \big( \| u(\cdot,s)\|^2_{V_\beta^2(\Omega)}
  + |s|^2  \| u(\cdot,s)\|^2_{V_\beta^0(\Omega)}\big)\, ds  \Big)^{1/2}
\]
are finite.  The proof of the following lemma is essentially the same as for nonweighted
spaces in \cite[Theorem 8.1]{av-64}. It is based on Plancherel's theorem for the Laplace transform
(see, e.~g.,  \cite[Formula (1.5.5)]{stenger}).

\begin{Le} \label{dl1}
Let $\gamma$ be a positive number. Then the Laplace transform realizes isomorphisms between the spaces
$L_{2,-\gamma}(0,\infty;V_\beta^0(\Omega))$ and $\stackrel{\circ}{W}\!{}_\beta^{l,1}(\Omega\times {\Bbb R_+},e^{-\gamma t})$ on one side
and $H_\beta(\Omega,\gamma)$ and $H_\beta^l(\Omega,\gamma)$, $l=1,2$, on the other side.
\end{Le}

\subsection{Solvability in \boldmath$\stackrel{\circ}{W}\!{}_\beta^{2,1}(\Omega\times {\Bbb R}_+,e^{-\gamma t})\times L_{2,-\gamma}(0,\infty;W_\beta^1(\Omega))$}

Using Theorem \ref{ct1}, we can easily prove the following theorem.

\begin{Th} \label{ct2}
Suppose that ${\mathfrak f}\in L_{2,-\gamma}(0,\infty;V_\beta^0(\Omega))$ and ${\mathfrak g}\in \stackrel{\circ}{W}\!{}_\beta^{1,1}(\Omega\times {\Bbb R}_+,e^{-\gamma t})$,
where the components $\beta_j$ of $\beta$ satisfy the inequalities
{\em (\ref{3ct1})} and $\gamma\ge\gamma_0$ with a sufficiently large positive number $\gamma_0$. Furthermore, we assume that ${\mathfrak g}$ satisfies the condition
\[
\int_\Omega {\mathfrak g}(x,t)\, dx = 0 \ \mbox{ for almost all }t.
\]
Then there exists a uniquely determined solution $({\mathfrak u},{\mathfrak p})$ of the problem {\em (\ref{stokes1})--(\ref{stokes3})} satisfying the estimate
\begin{eqnarray} \label{1ct2} \nonumber
&& \| {\mathfrak u}\|_{W_\beta^{2,1}(\Omega\times {\Bbb R}_+,e^{-\gamma t})} + \| {\mathfrak p}\|_{L_{2,-\gamma}({\Bbb R}_+;W_\beta^1(\Omega))}  \\
&& \le c\, \Big( \| {\mathfrak f}\|_{L_{2,-\gamma}({\Bbb R}_+;V_\beta^0(\Omega))} + \| {\mathfrak g} \|_{W_\beta^{1,1}(\Omega\times {\Bbb R}_+,e^{-\gamma t})} \Big)
\end{eqnarray}
and the condition
\begin{equation} \label{2ct2}
\int_\Omega {\mathfrak p}(x,t)\, dx =0 \ \mbox{for almost all }t.
\end{equation}
The constant $c$ is independent of $\gamma$ for $\gamma \ge\gamma_0$.
\end{Th}

Proof. Let $f \in H_\beta(\Omega,\gamma)$ and $g\in H_\beta^1(\Omega,\gamma)$ be the Laplace transforms of ${\mathfrak f}$ and ${\mathfrak g}$.
If $|s|\ge\mbox{Re}\, s > \gamma_0$, then there exists a unique solution
$(u,p)\in V_\beta^2(\Omega) \times (W_\beta^1(\Omega)\cap \stackrel{\circ}{L}\!\!{}_2(\Omega))$
of the problem (\ref{par3}), (\ref{par4}) satisfying the estimate (\ref{1ct1}) (see Theorem \ref{ct1}).
Integrating over the line $\mbox{Re}\, s = \sigma$
and taking the supremum with respect to $\sigma>\gamma$, we obtain the estimate (\ref{1ct2}) for the inverse Laplace transforms ${\mathfrak u},{\mathfrak p}$
of $u,p$. Obviously, the pair $({\mathfrak u},{\mathfrak p})$ is a solution of the problem (\ref{stokes1})--(\ref{stokes3}).
The uniqueness of the solution follows directly from Theorem \ref{ct1}. $\Box$

\subsection{Solvability in a finite \boldmath$t$-interval}

\begin{Th} \label{ct3}
Suppose that ${\mathfrak f}\in L_2(0,T;V_\beta^0(\Omega))$ and ${\mathfrak g}\in \stackrel{\circ}{W}\!{}_\beta^{1,1}(\Omega\times (0,T))$,
where the components $\beta_j$ of $\beta$ satisfy the inequalities {\em (\ref{3ct1})}. Furthermore, we assume that ${\mathfrak g}$ satisfies the condition
\[
\int_\Omega {\mathfrak g}(x,t)\, dx = 0 \ \mbox{ for almost all }t.
\]
Then there exists a uniquely determined solution $({\mathfrak u},{\mathfrak p})$ of the problem {\em (\ref{stokes1})--(\ref{stokes3})} satisfying the estimate
\begin{equation} \label{1ct3}
 \| {\mathfrak u}\|_{W_\beta^{2,1}(\Omega\times (0,T))} + \| {\mathfrak p}\|_{L_2(0,T;W_\beta^1(\Omega)}
 \le c\, \Big( \| {\mathfrak f}\|_{L_2(0,T;V_\beta^0(\Omega))}
  + \| {\mathfrak g} \|_{W_\beta^{1,1}(\Omega\times (0,T))} \Big)
\end{equation}
and the condition
\[
\int_\Omega {\mathfrak p}(x,t)\, dx =0 \ \mbox{for almost all }t.
\]
\end{Th}

Proof. Let ${\mathfrak F}\in L_2(0,\infty;V_\beta^0(\Omega))$ and ${\mathfrak G}\in \stackrel{\circ}{W}\!{}_\beta^{1,1}(\Omega\times {\Bbb R}_+)$
be extensions of ${\mathfrak f}$ and ${\mathfrak g}$ to $\Omega\times {\Bbb R}_+$ such that the integral of ${\mathfrak G}(\cdot,t)$
over $\Omega$ is zero for all $t$ and the estimates
\[
\| {\mathfrak F}\|_{L_2(0,\infty;V_\beta^0(\Omega))}\le c\, \| {\mathfrak f}\|_{L_2(0,T;V_\beta^0(\Omega))}\, ,\quad
\| {\mathfrak G}\|_{W_\beta^{1,1}(\Omega\times {\Bbb R}_+)} \le c\, \| {\mathfrak g}\|_{W_\beta^{1,1}(\Omega\times (0,T))}
\]
are satisfied. Obviously,  ${\mathfrak F}\in L_{2,-\gamma}(0,\infty;V_\beta^0(\Omega))$ and ${\mathfrak G}\in
\stackrel{\circ}{W}\!{}_\beta^{1,1}(\Omega\times {\Bbb R}_+,e^{-\gamma t})$
for arbitrary $\gamma \ge 0$. By Theorem \ref{ct2}, there exist a solution $({\mathfrak u},{\mathfrak p})$ of the problem (\ref{stokes1})--(\ref{stokes3}) satisfying the estimate (\ref{1ct2}) with ${\mathfrak F},{\mathfrak G}$ instead of ${\mathfrak f},{\mathfrak g}$ and arbitrary $\gamma \ge \gamma_0$. Then
\begin{eqnarray*}
&& \hspace{-2em}\| {\mathfrak u}\|_{W_\beta^{2,1}(\Omega\times (0,T))} + \|  {\mathfrak p}\|_{L_2(0,T;W_\beta^1(\Omega))}
\le e^{\gamma T} \, \Big( \| {\mathfrak u}\|_{W_\beta^{2,1}(\Omega\times {\Bbb R}_+,e^{-\gamma t})}
+ \|  {\mathfrak p}\|_{L_{2,-\gamma}(0,\infty;W_\beta^1(\Omega))}\Big) \\
&& \le c' e^{\gamma T}\, \Big( \| {\mathfrak F}\|_{L_{2,-\gamma}({\Bbb R}_+;V_\beta^0(\Omega))}
   + \| {\mathfrak G} \|_{W_\beta^{1,1}(\Omega\times {\Bbb R}_+,e^{-\gamma t})} \Big)\\
&&  \le c\, e^{\gamma T}\, \Big( \| {\mathfrak f}\|_{L_2(0,T;V_\beta^0(\Omega))}
  + \| {\mathfrak g} \|_{W_\beta^{1,1}(\Omega\times (0,T))} \Big),
\end{eqnarray*}
where $c$ is independent of $\gamma$. We show that the solution $({\mathfrak u},{\mathfrak p})$ depends only on ${\mathfrak f}$ and ${\mathfrak g}$ on $\Omega \times (0,T)$
but not on the extensions ${\mathfrak F}$ and ${\mathfrak G}$. Let $({\mathfrak F}',{\mathfrak G}')$ be another extension of $({\mathfrak f},{\mathfrak g})$, and let
$({\mathfrak u}',{\mathfrak p}')$ be the corresponding solution. Then Theorem \ref{ct2} yields
\begin{eqnarray*}
&& \| {\mathfrak u}-{\mathfrak u}'\|_{W_\beta^{2,1}(\Omega\times {\Bbb R}_+,e^{-\gamma t})}
+ \| {\mathfrak p}-{\mathfrak p}'\|_{L_{2,-\gamma}(0,\infty;W_\beta^1(\Omega))}\Big) \\
&& \le c \, \Big( \| {\mathfrak F}-{\mathfrak F}'\|_{L_{2,-\gamma}({\Bbb R}_+;V_\beta^0(\Omega))}
   + \| {\mathfrak G}-{\mathfrak G}' \|_{W_\beta^{1,1}(\Omega\times {\Bbb R}_+,e^{-\gamma t})} \Big)\\
&& \le c\, e^{-\gamma T}\, \Big( \| {\mathfrak F}-{\mathfrak F}'\|_{L_2(0,\infty;V_\beta^0(\Omega))}
   + \| {\mathfrak G}-{\mathfrak G}' \|_{W_\beta^{1,1}(\Omega\times {\Bbb R}_+)} \Big)
\end{eqnarray*}
for $\gamma > \gamma_0$, where $c$ is independent of $\gamma$. Since $e^{\gamma(T-\varepsilon-t)} \ge 1$ for $t\le T-\varepsilon$, we have
\[
\| {\mathfrak u}-{\mathfrak u}'\|_{W_\beta^{2,1}(\Omega\times (0,T-\varepsilon))}
  \le e^{\gamma(T-\varepsilon)}\, \| {\mathfrak u}-{\mathfrak u}'\|_{W_\beta^{2,1}(\Omega\times {\Bbb R}_+,e^{-\gamma t})}
\]
and
\[
\| {\mathfrak p}-{\mathfrak p}'\|_{L_2(0,T-\varepsilon;W_\beta^1(\Omega))}
  \le e^{\gamma(T-\varepsilon)} \, \| {\mathfrak p}-{\mathfrak p}'\|_{L_{2,-\gamma}(0,\infty;W_\beta^1(\Omega))}
\]
Consequently,
\begin{eqnarray*}
&& \| {\mathfrak u}-{\mathfrak u}'\|_{W_\beta^{2,1}(\Omega\times (0,T-\varepsilon))} + \| {\mathfrak p}-{\mathfrak p}'\|_{L_2(0,T-\varepsilon;W_\beta^1(\Omega))} \\
&& \le c\, e^{-\gamma\varepsilon}\, \Big( \| {\mathfrak F}-{\mathfrak F}'\|_{L_2(0,\infty;V_\beta^0(\Omega))}
   + \| {\mathfrak G}-{\mathfrak G}' \|_{W_\beta^{1,1}(\Omega\times {\Bbb R}_+)} \Big),
\end{eqnarray*}
where $c$ is independent of $\gamma$ for $\gamma> \gamma_0$. If we let $\gamma$ tend to infinity, we conclude that ${\mathfrak u}={\mathfrak u}'$
and ${\mathfrak p}={\mathfrak p}'$ in $\Omega\times (0,T-\varepsilon)$ for arbitrary $\varepsilon >0$.  This proves the theorem. $\Box$

\end{document}